\documentclass{llncs}
\usepackage{llncsdoc}
\usepackage{amsmath}
\usepackage{amssymb}
\usepackage{footnote}
\newcommand{\FF}{\mathbb{F}}
\begin{document}
\title{Efficient Methods in Counting Generalized Necklaces
%\thanks{The author would like to thank Professor Aaron Gulliver,
%Department of Electrical and Computer Engineering,
%University of Victoria, Canada for introducing him to this problem.}\\
\hfill \\}
\author{V. Ch. Venkaiah} \institute{School of Computer and
Information Sciences\\ University of Hyderabad, Gachibowli
\\ Hyderabad 500 046, India, \\ \email{venkaiah@hotmail.com;
vvcs@uohyd.ernet.in}\\ \vspace{0.1in}}

\maketitle

%\hspace{-0.20in}{\sl Dedicated to my wife Mrs. Anuradha, who has
%been with me throughout.}
\begin{abstract}
It is shown in \cite{venk14} by Venkaiah in 2015 that a category
of the number of generalized necklaces can be computed using the
expression
\begin{equation*}
e(n, q) = \frac{1}{(q-1) ord(\lambda)
n}\sum^{ord(\lambda)n}_{\substack{t \in \FF_{q} \setminus \{0\},
i=1
\\ t^{\frac{n}{\gcd(n, i)}} \lambda^{\frac{i}{\gcd(n,i)}} =
1}}(q^{\gcd(n,i)} - 1) + 1
\end{equation*}
where $q$ (number of colors) is the size of the prime field
$\FF_{q}$, $\lambda$ is the constant of the consta cyclic shift,
$n$ is the length of the necklace. However, direct evaluation of
this expression requires, apart from the $\gcd$ computations,
$2*(q-1)*Ord(\lambda)*n$ exponentiations and
$(q-1)*Ord(\lambda)*n$ multiplications, at most
$(q-1)*Ord(\lambda)*n$ exponentiations and at most
$2*(q-1)*Ord(\lambda)*n$ additions and hence computationally
intensive. This note discusses various other ways of evaluating
the expression and tries to throw some light on amortizing the
amount of computation.
\end{abstract}
\section{Introduction}
A $q-$ary necklace of length $n$ is an equivalence class of
$q-$ary strings of length $n$ under rotation. Various other forms
of necklaces exist in the literature. Unlabelled necklaces and
Lyndon words are some such. The number of $q-$ary necklaces of
length $n$ is given by the formula \cite{berstel07,bou10,fred77}
\begin{equation*}
N_{q}(n) = \frac{1}{n} \sum_{d \mid n}\varphi(m/d)q^{d}
\end{equation*}
where $\varphi$ is the Euler's totient function.

Three categories of generalizations to the necklaces are proposed
in \cite{venk14}. First category of generalized necklace is an
equivalence class of $q-$ary strings under rotation as well as
multiplication by a nonzero constant. Second category of
generalized necklaces is an equivalence class of $q-$ary strings
under consta cyclic shift with a fixed nonzero element $\lambda
\in \FF_{q}\setminus\{0\}$. The third category of generalized
necklaces is an equivalence class of $q-ary$ strings under consta
cyclic shift with some nonzero constant $\lambda$ as well as
multiplication by a nonzero constant. Following are examples of
these generalized necklaces.

\hspace{-0.25in}{\bf Example 1} Let $q = 5$ and $n = 2$.
Representing $5-$ary strings of length 2 as the first degree
polynomials over $\FF_{5}$, we have the set $\{0, 1, 2, 3, 4, x,
x+1, x+2, x+3, x+4, 2x, 2x+1, 2x+2, 2x+3, 2x+4, 3x, 3x+1, 3x+2,
3x+3, 3x+4, 4x, 4x+1, 4x+2, 4x+3, 4x+4\}$. Define a relation that
two elements of this set are related if one is a cyclic shift of
the other or it is a non-zero constant multiple of the other or
both. This gives rise to the following five equivalence classes:
$E_{1} = \{1, x, 2, 2x, 4, 4x, 3, 3x\}$, $E_{2} =\{x+1, 2x+2,
3x+3, 4x+4\}$, $E_{3} = \{x+2, 2x+1, 2x+4, 3x+1, 4x+3, 4x+2, x+3,
3x+4\}$, and $E_{4} = \{x+4, 4x+1, 2x+3, 3x+2\}$, $E_{5} = \{ 0
\}$. Representatives of each of these equivalence classes will be
the first category of generalized necklaces.

\hspace{-0.25in}{\bf Example 2} Let $q = 11$ and $n = 2$.
Representing $11-$ary strings of length 2 as the first degree
polynomials over $\FF_{11}$, we have the set  $ \{0, 1, 2, 3, 4,
5, \\ 6, 7, 8, 9, 10, x, x+1, x+2, x+3, x+4, x+5, x+6, x+7, x+8,
x+9, x+10, 2x, 2x+1, 2x+2, 2x+3, 2x+4, 2x+5, 2x+6, 2x+7, 2x+8,
2x+9, 2x+10, 3x, 3x+1, 3x+2, 3x+3, 3x+4, 3x+5, 3x+6, 3x+7, 3x+8,
3x+9, 3x+10, 4x, 4x+1, 4x+2, 4x+3, 4x+4, 4x+5, 4x+6, 4x+7, 4x+8,
4x+9, 4x+10, 5x, 5x+1, 5x+2, 5x+3, 5x+4, 5x+5, 5x+6, 5x+7, 5x+8,
5x+9, 5x+10, 6x, 6x+1, 6x+2, 6x+3, 6x+4, 6x+5, 6x+6, 6x+7, 6x+8,
6x+9, 6x+10, 7x, 7x+1, 7x+2, 7x+3, 7x+4, 7x+5, 7x+6, 7x+7, 7x+8,
7x+9, 7x+10, 8x, 8x+1, 8x+2, 8x+3, 8x+4, 8x+5, 8x+6, 8x+7, 8x+8,
8x+9, 8x+10, 9x, 9x+1, 9x+2, 9x+3, 9x+4, 9x+5, 9x+6, 9x+7, 9x+8,
9x+9, 9x+10, 10x, 10x+1, 10x+2, 10x+3, 10x+4, 10x+5, 10x+6, 10x+7,
10x+8, 10x+9, 10x+10\} $. Define a relation that two elements of
this set are related if one is a cyclic shift of the other or it
is a non-zero constant multiple of the other or both. This gives
rise to the following eight equivalence classes: $E_{1} = \{1, x,
2, 2x, 3, 3x, 4, 4x, 5, 5x, 6, 6x, 7, 7x, 8, 8x, 9, 9x, 10,
10x\}$, $E_{2} = \{ x+1, 2x+2, 3x+3, 4x+4, 5x+5, 6x+6, 7x+7, 8x+8,
9x+9, 10x + 10\}$, $E_{3} = \{x+2, 2x+1, 2x+4, 3x+6, 4x+8, 5x+10,
6x+1, 7x+3, 8x+5, 9x+7, 10x+9, 4x+2, 6x+3, 8x+4, 10x+5, x+6, 3x+7,
5x+8, 7x+9, 9x+10\}$, $E_{4} = \{x+3, 3x+1, 2x+6, 3x+9, 4x+1,
5x+4, 6x+7, 7x+10, 8x+2, 9x+5, 10x+8, 6x+2, 9x+3, x+4, 4x+5, 7x+6,
10x+7, 2x+8, 5x+9, 8x+10\}$, $E_{5} = \{x+5, 5x+1, 2x+10, 3x+4,
4x+9, 5x+3, 6x+8, 7x+2, 8x+7, 9x+1, 10x+6, 10x+2, 4x+3, 9x+4,
3x+5, 8x+6, 2x+7, 7x+8, x+9, 6x+10\}$, $E_{6} = \{x+7, 7x+1, 2x+3,
3x+10, 4x+6, 5x+2, 6x+9, 7x+5, 8x+1, 9x+8, 10x+4, 3x+2, 10x+3,
6x+4, 2x+5, 9x+6, 5x+7, x+8, 8x+9, 4x+10\}$, $E_{7} = \{x+10,
10x+1, 2x+9, 3x+8, 4x+7, 5x+6, 6x+5, 7x+4, 8x+3, 9x+2\}$, $E_{8} =
\{0\}$. As in the previous example, representatives of each of
these equivalence classes will be the first category of
generalized necklaces.

\hspace{-0.25in}{\bf Example 3} Let $q = 5$, $n = 2$, and $\lambda
= 2$. Representing the elements as the first degree polynomials we
have the set given in Example 1. Define a relation that two
elements are related if one is a consta cyclic shift with $\lambda
= 2$. This gives rise to the following four equivalence classes of
$S_{2}$: $E_{1} = \{1, x, 2, 2x, 4, 4x, 3, 3x\}$, $E_{2} =\{x+1,
x+2, 2x+2, 2x+4, 4x+4, 4x+3, 3x+3, 3x+1\}$, $E_{3} = \{x+3, 3x+2,
2x+1, x+4, 4x+2, 2x+3, 3x+4, 4x+1\}$, and $E_{4} = \{ 0 \}$.
Representatives of each of these equivalence classes will be the
second category of generalized necklaces.

\hspace{-0.25in}{\bf Example 4} Let $q = 11$, $n = 2$, and
$\lambda = 3$. Representing the elements as the first degree
polynomials, we have the set given in Example 2. Defining a
relation as in Example 3, we have the following fourteen
equivalence classes: $E_{1} = \{1, x, 3, 3x, 9, 9x, 5, 5x, 4,
4x\}$, $E_{2} = \{2, 2x, 6, 6x, 7, 7x, 10, 10x, 8, 8x\}$, $E_{3} =
\{x+1, x+3, 3x+3, 3x+9, 9x+9, 9x+5, 5x+5, 5x+4, 4x+4, 4x+1\}$,
$E_{4} = \{x+2, 2x+3, 3x+6, 6x+9, 9x+7, 7x+5, 5x+10, 10x+4, 4x+8,
8x+1\}$, $E_{5} = \{x+4, 4x+3, 3x+1, x+9, 9x+3, 3x+5, 5x+9, 9x+4,
4x+5, 5x+1\}$, $E_{6} = \{x+5, 5x+3, 3x+4, 4x+9, 9x+1\}$, $E_{7} =
\{x+6, 6x+3, 3x+7, 7x+9, 9x+10, 10x+5, 5x+8, 8x+4, 4x+2, 2x+1\}$,
$E_{8} = \{x+7, 7x+3, 3x+10, 10x+9, 9x+8, 8x+5, 5x+2, 2x+4, 4x+6,
6x+1\}$, $E_{9} = \{x+8, 8x+3, 3x+2, 2x+9, 9x+6, 6x+5, 5x+7, 7x+4,
4x+10, 10x+1\}$, $E_{10} = \{x+10, 10x+3, 3x+8, 8x+9, 9x+2, 2x+5,
5x+6, 6x+4, 4x+7, 7x+1\}$, $E_{11} = \{2x+2, 2x+6, 6x+6, 6x+7,
7x+7, 7x+10, 10x+10, 10x+8, 8x+8, 8x+2\}$, $E_{12} = \{2x+7, 7x+6,
6x+10, 10x+7, 7x+8, 8x+10, 10x+2, 2x+8, 8x+6, 6x+2\}$, $E_{13} =
\{2x+10, 10x+6, 6x+8, 8x+7, 7x+2\}$, $E_{14} = \{0\}$.
Representatives of each of these equivalence classes will be the
second category of generalized necklaces.

\hspace{-0.25in}{\bf Example 5} Let $q$, $n$, and $\lambda$ be as
in the Example 4, so that the set of elements is as in that
example. Define a relation that two elements of this set are
related if one is a consta cyclic shift with $\lambda = 3$ of the
other or it is a non-zero multiple of the other or both. This
gives rise to the following eight equivalence classes: $E_{1} =
\{1, x, 3, 3x, 9, 9x, 5, 5x, 4, 4x, 2, 2x, 6, 6x, 7, 7x, 10, 10x,
8, 8x\}$, $E_{2} = \{x+1, x+3, 3x+3, 3x+9, 9x+9, 9x+5, 5x+5, 5x+4,
4x+4, 4x+1, 2x+2, 2x+6, 6x+6, 6x+7, 7x+7, 7x+10, 10x+10, 10x+8,
8x+8, 8x+2\}$, $E_{3} = \{x+2, 2x+3, 3x+6, 6x+9, 9x+7, 7x+5,
5x+10, 10x+4, 4x+8, 8x+1, x+7, 7x+3, 3x+10, 10x+9, 9x+8, 8x+5,
5x+2, 2x+4, 4x+6, 6x+1\}$, $E_{4} = \{x+4, 4x+3, 3x+1, x+9, 9x+3,
3x+5, 5x+9, 9x+4, 4x+5, 5x+1, 2x+7, 7x+6, 6x+10, 10x+7, 7x+8,
8x+10, 10x+2, 2x+8, 8x+6, 6x+2\}$, $E_{5} = \{x+5, 5x+3, 3x+4,
4x+9, 9x+1, 2x+10, 10x+6, 6x+8, 8x+7, 7x+2\}$, $E_{6} = \{x+6,
6x+3, 3x+7, 7x+9, 9x+10, 10x+5, 5x+8, 8x+4, 4x+2, 2x+1\}$, $E_{7}
= \{x+8, 8x+3, 3x+2, 2x+9, 9x+6, 6x+5, 5x+7, 7x+4, 4x+10, 10x+1,
x+10, 10x+3, 3x+8, 8x+9, 9x+2, 2x+5, 5x+6, 6x+4, 4x+7, 7x+1\}$,
$E_{8} = \{0\}$. Representatives of each of these equivalence
classes will be the third category of generalized necklaces.

The following expression given in \cite{venk13,venk14}
\begin{equation}
b(n, q) = \frac{1}{(q-1)n}\sum_{d \mid n} \varphi(d) \gcd(d, q-1)
(q^{n/d} - 1) +1,
\end{equation}
and the expression given in \cite{venk14}
\begin{equation*}
d(n, q) =
\frac{1}{ord(\lambda)n}\sum^{ord(\lambda)n}_{\substack{i=1
\\ \lambda^{\frac{i}{\gcd(n,i)}} =
1}}(q^{\gcd(n,i)} - 1) + 1,
\end{equation*}
can be used to count the number of first and second category of
generalized necklaces respectively.

While it is computationally easy to evaluate these expressions
given for counting the number of necklaces, first, and second
category of generalized necklaces, it can be seen from the example
given in the next section that it is computationally demanding to
evaluate the expression that counts the number of third category
of generalized necklaces. The purpose of this paper is to address
this problem and propose alternative ways of evaluating the
expression and suggest techniques that amortize the amount of
computation.
\section{Third Category Generalized Necklaces}
As mentioned previously, it is an equivalence class of $q-$ary
strings under consta cyclic shift as well as multiplication by a
nonzero constant. This class of necklaces and their number come up
in the study of quasi twisted codes \cite{gulliver}. A quasi
twisted code of index $p$ is a code in which a consta cyclic shift
of a codeword by $p$ places is another codeword. It may be noted
that a quasi twisted code generalizes the class of consta cyclic
codes ($p=1$), quasi-cyclic codes ($\lambda = 1$) and the class of
cyclic codes ($\lambda =1, p = 1$.)\cite{Chen08}.

The following expression given in \cite{venk14} can be used to
count the number of these generalized necklaces.
\begin{equation*}
e(n, q) = \frac{1}{(q-1) ord(\lambda)
n}\sum^{ord(\lambda)n}_{\substack{t \in \FF_{q} \setminus \{0\},
i=1
\\ t^{\frac{n}{\gcd(n, i)}} \lambda^{\frac{i}{\gcd(n,i)}} =
1}}(q^{\gcd(n,i)} - 1) + 1
\end{equation*}
Following example illustrates a way of evaluating the expression
to count these necklaces:

\hspace{-0.25in}{\bf Example 6} In Example 4, $q = 11$, $n = 2$,
and $\lambda = 3$. Since $5$ is the least integer so that,
$\lambda^{5} = 3^{5} = 1$, we have $Ord(\lambda) = Ord(3) = 5$.
Therefore,
\begin{eqnarray*}
e(2, 11) &=& \frac{1}{100} \sum_{\substack{t \in \FF_{11} \setminus \{0\}, i = 1 \\
t^{\frac{2}{\gcd(2, i)}} 3^{\frac{i}{\gcd(2, i)}} =
1}}^{10}(11^{\gcd(2, i)} - 1) + 1
 \nonumber \\
 &=& \frac{1}{100}\left\{ \left[ (11^{\gcd(2, 5)} - 1) +
(11^{\gcd(2, 10)} - 1)\right] + \left[ (11^{\gcd(2, 1)} -
1)\right] \right\}
 \nonumber \\
 \end{eqnarray*}
 \begin{eqnarray*}
 \hspace{1cm}
&&  + \frac{1}{100} \left\{\left[ (11^{\gcd(2,3)} - 1) +
(11^{\gcd(2, 8)} - 1) \right] +
\left[ (11^{\gcd(2, 2)} - 1) + (11^{\gcd(2, 7)} - 1) \right] \right\} \nonumber \\
&&  + \frac{1}{100} \left\{\left[ (11^{\gcd(2, 4)} - 1) +
(11^{\gcd(2, 9)} - 1) \right] +
\left[ (11^{\gcd(2, 9)} - 1) \right] \right\} \nonumber \\
&& + \frac{1}{100} \left\{ \left[ (11^{\gcd(2, 7)} - 1) \right] +
\left[ (11^{\gcd(2,
3)} - 1) \right] \right\} \nonumber \\
&&  + \frac{1}{100} \left\{ \left[ (11^{\gcd(2, 1)} - 1) +
(11^{\gcd(2, 6)} - 1) \right] + \left[ (11^{\gcd(2, 5)} - 1)
\right] \right\} + 1.
\end{eqnarray*}
This is because
\begin{eqnarray*}
1^{\frac{2}{\gcd(2,5)}}3^{\frac{5}{\gcd(2, 5)}} =
1^{\frac{2}{\gcd(2,10)}}3^{\frac{10}{\gcd(2, 10)}} = \\
2^{\frac{2}{\gcd(2,1)}}3^{\frac{1}{\gcd(2, 1)}} = \\
3^{\frac{2}{\gcd(2,3)}}3^{\frac{3}{\gcd(2, 3)}} =
3^{\frac{2}{\gcd(2,8)}}3^{\frac{8}{\gcd(2, 8)}} = \\
4^{\frac{2}{\gcd(2,2)}}3^{\frac{2}{\gcd(2, 2)}} =
4^{\frac{2}{\gcd(2,7)}}3^{\frac{7}{\gcd(2, 7)}} = \\
5^{\frac{2}{\gcd(2,4)}}3^{\frac{4}{\gcd(2, 4)}} =
5^{\frac{2}{\gcd(2,9)}}3^{\frac{9}{\gcd(2, 9)}} = \\
6^{\frac{2}{\gcd(2,9)}}3^{\frac{9}{\gcd(2, 9)}} = \\
7^{\frac{2}{\gcd(2,7)}}3^{\frac{7}{\gcd(2, 7)}} = \\
8^{\frac{2}{\gcd(2,3)}}3^{\frac{3}{\gcd(2, 3)}} = \\
9^{\frac{2}{\gcd(2,1)}}3^{\frac{1}{\gcd(2, 1)}} =
9^{\frac{2}{\gcd(2,6)}}3^{\frac{6}{\gcd(2, 6)}} = \\
10^{\frac{2}{\gcd(2,5)}}3^{\frac{5}{\gcd(2, 5)}} = 1
\end{eqnarray*}
So,
\begin{equation*}
e(2, 11) = \frac{1}{100} [ 130 + 10 + 130 + 130 + 130 + 10 + 10 +
10 + 130 + 10 ] + 1 = 8.
\end{equation*}

Note that the evaluation required the computation of
$t^{\frac{2}{\gcd(2,i)}}3^{\frac{i}{\gcd(2,i)}}$ for every $t \in
\FF_{11} \setminus \{0\}$ and $1 \leq i \leq 10$. So, the
computational effort required to evaluate the expression, apart
from the $\gcd$ computations, is 2*10*10 exponentiations  and
10*10 multiplications to compute $t^{\frac{2}{\gcd(2, i)}}
3^{\frac{i}{\gcd(2, i)}}$ for $t \in \FF_{11} \setminus \{0\}$ and
$1 \leq i \leq 10$, at most 10*10 exponentiations and at most
2*10*10 additions to compute $11^{\gcd(2, i)}$ and their sum for
$t$ and $i$ such that
$t^{\frac{2}{\gcd(2,i)}}3^{\frac{i}{\gcd(2,i)}} = 1$, where $t \in
\FF_{11} \setminus \{0\}$ and $1 \leq i \leq 10$, and couple of
multiplications, divisions, and additions or subtractions.
Generalizing these, the computational complexity to evaluate the
expression can be seen to require $2*(q-1)*Ord(\lambda)*n$
exponentiations and $(q-1)*Ord(\lambda)*n$ number of
multiplications to compute $t^{\frac{n}{\gcd(n, i)}}
\lambda^{\frac{i}{\gcd(n, i)}}$ for $t \in \FF_{q} \setminus
\{0\}$ and $1 \leq i \leq Ord(\lambda)*n$, at most
$(q-1)*Ord(\lambda)*n$ exponentiations and at most
$2*(q-1)*Ord(\lambda)*n$ additions to compute $q^{\gcd(n, i)}$ and
their sum for $t$ and $i$ such that
$t^{\frac{n}{\gcd(n,i)}}\lambda^{\frac{i}{\gcd(n,i)}} = 1$, where
$t \in \FF_{q} \setminus \{0\}$ and $1 \leq i \leq
Ord(\lambda)*n$, and couple of multiplications, divisions, and
additions or subtractions. This indicates that the effort required
is quite huge even for moderate field sizes.
\section{Proposed Method}
The following algorithm is an improvement over the foregoing
method in that it computes $\lambda^{\frac{i}{\gcd(n, i)}}$ for
every $i$ only once as against computing them afresh for each
distinct $t$.
\subsection{Algorithm}
\begin{itemize}
\item Step 1: Compute $\lambda-exp[i] = \lambda^{\frac{i}{\gcd(n,
i)}}$ for $1 \leq i \leq Ord(\lambda)n$ and store. \item Step 2:
Let $count = 0$. \item Step 3: For a fixed $i$ and for a fixed $t$
\item \hspace{0.4cm} Step 3.1: Compute $t_{exp} =
t^{\frac{n}{\gcd(n,i)}}$  \item \hspace{0.4cm} Step 3.2 Compute
$t\lambda-exp[i] = t_{exp}*\lambda-exp[i]$ by fetching $\lambda
-exp[i]$ from the memory. \item \hspace{0.4cm} Step 3.3: If
$t\lambda-exp[i] = 1$ then compute $count = count + q^{\gcd(n,i)}
- 1$. \item Step 4: Repeat Step 3 for every $i, 1 \leq i \leq
Ord(\lambda)n$ and for every $t \in \FF_{q} \setminus\{0\}$. \item
Step 5: Compute $count = (count / ((q-1)ord(\lambda)n)) + 1$.
\end{itemize}
Computational effort required by this method, apart from the
$\gcd$ computations, can be seen to be the following:
$Ord(\lambda)n$ exponentiations for Step 1, $((q-1)Ord(\lambda)n)$
exponentiations, $((q-1)Ord(\lambda)n)$ multiplications, at most
$((q-1)Ord(\lambda)n)$ exponentiations, and at most
$2*((q-1)Ord(\lambda)n)$ additions or subtractions for Steps 3 and
4, and couple of multiplications, divisions, and an addition for
Step 5.

Note that there is approximately $50\% $ reduction in the number
of exponentiations. However, even this can be seen to be quite
huge. So, we propose the following method that exploits the
multiplicity of the $\gcd$s to bring down the computational
requirement.
\subsection{Algorithm}
\begin{itemize}
\item Step 1: for $1 \leq i \leq Ord(\lambda)n$ \\
\hspace{   0.5cm} Compute $\gcd(n, i)$ and store at the $i^{th}$
location of the array $gcdiv$. That is $gcdiv[i] = \gcd(n, i)$.\\
\item Step 2: Let $g_{i}, 1 \leq i \leq d,$ be the distinct of
these gcds. That is, $g_{j} = \gcd(n, k)$ for some $k, 1 \leq k
\leq Ord(\lambda)n$, $g_{i} \neq g_{j}$ for $i \neq j$, and for
any $k, 1 \leq k \leq Ord(\lambda)n$ there is an
$\ell, 1 \leq \ell \leq d,$ such that $\gcd(n, k) = g_{\ell}$. \\
\item Step 3: Also, let $G_{j} = \{i : 1 \leq i \leq Ord(\lambda)n
\mbox{ and } \gcd(n, i) = g_{j}\}$ for $1 \leq j \leq d$. \\ \item
Step 4: for $1 \leq i \leq Ord(\lambda)n$
\\ Compute
$\lambda^{\frac{i}{\gcd(n, i)}}$ and store at the $i^{th}$
location of the array $\lambda-exp$. That is, $\lambda-exp[i] =
\lambda^{\frac{i}{\gcd(n,i)}}$. \\ \item Step 5: for $1 \leq i
\leq d$ \\ Compute $\frac{n}{g_{i}}$ and store at the $i^{th}$
location of the array $nbygcd$. That is, $nbygcd[i] =
\frac{n}{g_{i}}$. \\ \item Step 6: $count = 0.0$
\\ \mbox{         } for
$1 \leq t < q$ \\\mbox{                  } for $1 \leq j \leq d$ \\
\mbox{                           } Compute $texp[j] =
t^{nbygcd[j]}$ \\ \mbox{                           } for every $k
\in G_{j}$ \\ \mbox{                           } begin \\ \mbox{
      } Compute $t\lambda-exp$ = $texp[j]*\lambda-exp[k]$
\\\mbox{                                       } if $t\lambda-exp =
1$ then $count = count + (q^{g_{j}} - 1)$. \\ \mbox{
             } end \item Step 7: Number of Third Category generalized Necklaces
= $(count / ((q-1)Ord(\lambda)n)) + 1$
\end{itemize}

\hspace{-0.5cm} \emph{\bf Example 7}
\begin{itemize}
\item $n = 2$, $\lambda = 3$, $q = 11$, so $Ord(\lambda) = Ord(3)
= 5$, and $Ord(\lambda)n = 10$\\
\item $\gcd(2, 1) = 1$, $\gcd(2, 2) = 2$, $\gcd(2, 3) = 1$,
$\gcd(2, 4) = 2$, $\gcd(2, 5) = 1$, $\gcd(2, 6) = 2$, $\gcd(2, 7)
= 1$, $\gcd(2, 8) = 2$, $\gcd(2, 9) = 1$, $\gcd(2, 10) = 2$ \\
\item $g_{1} = 1$, $g_{2} = 2$ \\
\item $G_{1} = \{1, 3, 5, 7, 9\}$ and $G_{2} = \{2, 4, 6, 8, 10\}$
\\
\item $\lambda-exp[1] = 3^{\frac{1}{1}} = 3$; $\lambda-exp[2] =
3^{\frac{2}{2}} = 3$; $\lambda-exp[3] = 3^{\frac{3}{1}} = 5$;
$\lambda-exp[4] = 3^{\frac{4}{2}} = 9$; $\lambda-exp[5] =
3^{\frac{5}{1}} = 1$; $\lambda-exp[6] = 3^{\frac{6}{2}} = 5$;
$\lambda-exp[7] = 3^{\frac{7}{1}} = 9$; $\lambda-exp[8] =
3^{\frac{8}{2}}
= 4$; $\lambda-exp[9] = 3^{\frac{9}{1}} = 4$; $\lambda-exp[10] = 3^{\frac{10}{2}} = 1$ \\
\item $nbygcd[1] = 2/1 = 2$; $nbygcd[2] = 2/2 = 1$ \\
\item $t = 1$; $j = 1$; $texp[1] = t^{nbygcd[1]} = t^{2} = 1$;
$t\lambda-exp = 1$ for $k = 5$; So, $count = (q^{g_{1}} -
1) = (11^{1} - 1) = 10$ \\
\item $t = 1$; $j = 2$; $texp[2] = t^{nbygcd[2]} = t^{1} = 1$;
$t\lambda-exp = 1$ for $k = 10$; So, $count = count + (q^{g_{2}} -
1) = 10 + (11^{2} - 1) = 130$ \\
\item $t = 2$; $j = 1$; $texp[1] = t^{nbygcd[1]} = t^{2} = 4$;
$t\lambda-exp = 1$ for $k = 1$; So, $count = count + (q^{g_{1}} -
1) = 130 + (11^{1} - 1) = 140$ \\
\item $t = 2$; $j = 2$; $texp[2] = t^{nbygcd[2]} = t^{1} = 2$;
$t\lambda-exp \neq 1$ for any $k \in G_{2}$; So, count remains same as above. \\
\item $t = 3$; $j = 1$; $texp[1] = t^{nbygcd[1]} = t^{2} = 9$;
$t\lambda-exp = 1$ for $k = 3$; So, $count = count + (q^{g_{1}} -
1) = 140 + (11^{1} - 1) = 150$ \\
\item Continuing in this manner, we have \\
\item $t = 10$; $j = 1$; $texp[1] = t^{nbygcd[1]} = t^{2} = 1$;
$t\lambda-exp = 1$ for $k = 5$; So, $count = count + (q^{g_{1}} -
1) = 690 + (11^{1} - 1) = 700$ \\
\item $t = 10$; $j = 2$; $texp[2] = t^{nbygcd[2]} = t^{1} = 10$;
$t\lambda-exp \neq 1$ for any $k \in G_{2}$; So count remains same as above.\\
\item Number of third category of generalized necklaces is
$\frac{count}{(q-1)Ord(\lambda)n} + 1 = \frac{700}{10*5*2} + 1 =
8$
\end{itemize}

Computational effort required by the proposed method given in Sec.
3.2 is as follows: Step 1 requires $Ord(\lambda)n$ $\gcd$
computations, Steps 2 and 3 require at most $d*Ord(\lambda)n$
operations, where $d$ is the number of distinct $\gcd$s. Step 4
requires $Ord(\lambda)n$ exponentiations, Step 5 requires $d$
divisions, Step 6 requires $(q-1)d$ exponentiations,
$(q-1)Ord(\lambda)n$ multiplications, at most $(q-1)Ord(\lambda)n$
exponentiations, and at most $2*(q-1)Ord(\lambda)n$ additions or
subtractions. Step 7 requires couple of divisions and additions.
\section{Conclusions}
Couple of methods of evaluating the expression to count the third
category of generalized necklaces are studied; while doing so, we
formalized the ideas into algorithms. Comparative analysis of
their computational complexity is carried out and found that the
best alternative method suggested in the paper cuts down on the
amount of computation by more than $50\%$. Research is underway to
arrive at still better efficient methods of evaluating the
expression.

\end{document}